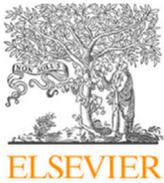

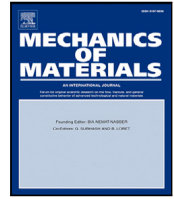

Editor invited article

# An efficient active-stress electromechanical isogeometric shell model for muscular thin film simulations

Michele Torre [a], Simone Morganti [b], Alessandro Nitti [c], Marco Donato de Tullio [c], Josef Kiendl [d], Francesco Silvio Pasqualini [a], Alessandro Reali [a,*]

[a] *Department of Civil Engineering and Architecture, University of Pavia, via Ferrata 3, Pavia, 27100, Italy*
[b] *Department of Electrical, Computer, and Biomedical Engineering, University of Pavia, via Ferrata 5, Pavia, 27100, Italy*
[c] *Department of Mechanics, Mathematics and Management, Polytechnic University of Bari, Via Re David 200, Bari, 70125, Italy*
[d] *Institute of Engineering Mechanics and Structural Analysis, Bundeswehr University Munich, Werner-Heisenberg-Weg 39, Neubiberg, 85577, Germany*

## ARTICLE INFO

## ABSTRACT

We propose an isogeometric approach to model the deformation of active thin films using layered, nonlinear, Kirchhoff–Love shells. Isogeometric Collocation and Galerkin formulations are employed to discretize the electrophysiological and mechanical sub-problems, respectively, with the possibility to adopt different element and time-step sizes.

Numerical tests illustrate the capabilities of the active-stress-based approach to effectively simulate the contraction of thin films in both quasi-static and dynamic conditions.

## 1. Introduction

The quantification of tissue-generated forces during contraction is of utmost relevance in many biomedical applications (Schroer et al., 2017; MacQueen et al., 2018). However, direct measurements are extremely complex and, therefore, researchers often resort to assessing the deformations of artificial constructs composed of a biological tissue and a substrate with known geometry and material properties (Feinberg et al., 2007; Grosberg et al., 2011), which are successively correlated with force generation via different mechanistic models. For instance, the curvature of a muscular thin film (MTF), measured by sensors embedded in the construct (Lind et al., 2017), can be related to the state of stress via either the linear modified Stoney equation (Feinberg et al., 2007), other non-linear models (Alford et al., 2010) or finite element simulations (Shim et al., 2012; Pezzuto et al., 2014).

MTFs are flat membranes composed of a few layers of cardiomyocytes grown on top of an elastic substrate, which serves as a supporting structure, favoring the cell alignment and organization in fibers (Agarwal et al., 2013), which depends on the topography of the cell-to-substrate interface. In this context, standard simulations employ a complete three-dimensional description, even if the mechanical response of the constructs has a low-order dimensionality due to the high aspect ratio of the domain shape, as highlighted for beam-like structures (Nardinocchi et al., 2015), which can be effectively discretized via isogeometric analysis (Ferri et al., 2023). We believe that a structural shell model can effectively reproduce the coupled electromechanical behavior of thin structures, reducing the computational effort to the simulation of a bi-variate manifold rather than a three-dimensional problem. Furthermore, a similar rationale holds for the electrophysiological part of the problem: Since the biological tissue is thin, the monodomain formulation used to simulate the cell activity can be solved on a bi-variate manifold.

Herein, we develop an alternative approach to finite elements to simulate tissue contractility focusing on MTFs of the type presented in Shim et al. (2012), as they are simple and relevant case studies, considering also simulations of similar constructs (Lucantonio et al., 2014) – for instance, those with curved geometry (MacQueen et al., 2018; Böl et al., 2009) – can benefit from the same development.

Focusing on Kirchhoff–Love shells, we introduce the active stress formulation to simulate cell activation. We provide the dual alternative to the active strain approach, which is used to model tissue contractility in shells (Nitti et al., 2021), for the cases where the constitutive model is already calibrated according to the former approach (Shim et al., 2012). Moreover, differently from previous works, we use two different numerical approaches to discretize the electrophysiological and mechanical sub-problems combining an Isogeometric Collocation (Auricchio et al., 2012) and an Isogeometric Galerkin (Kiendl et al., 2009)

---

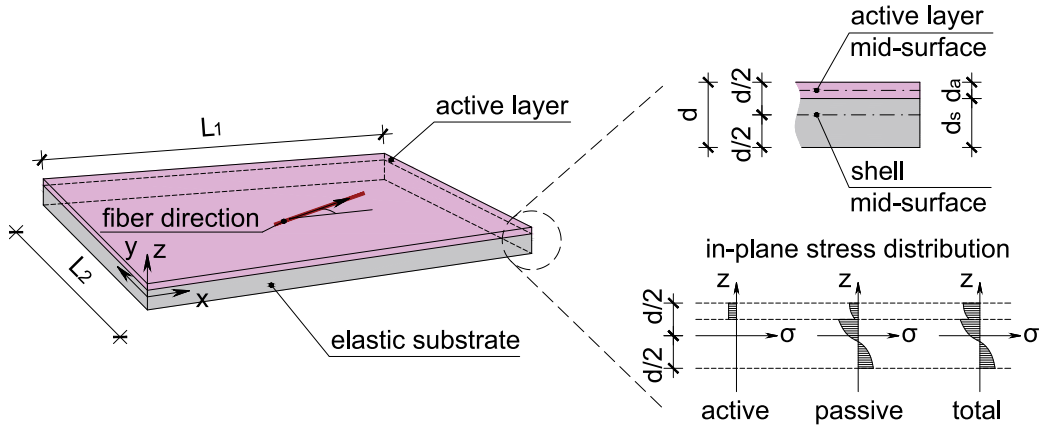

**Fig. 1.** Schematic representation of MTFs composed by an active biological layer and an elastic substrate, both represented in the numerical approach by their mid-surfaces. Since the two materials have different mechanical properties, a discontinuity in the in-plane stresses is generated by the compatibility of displacements.

formulation, taking advantage of B-spline continuity (Cottrell et al., 2009). On the one hand, in electrophysiology, the continuity enables a collocation approach that limits the computational effort due to the computation of the reactive term of the monodomain formulation, while preserving the accuracy (Torre et al., 2022). On the other hand, in mechanics, the continuity allows to efficiently discretize curvatures needed for the Kirchhoff–Love shell model (Kiendl et al., 2009). Altogether, the resulting approach guarantees a high accuracy at a relatively low cost also in geometrically-complex situations.

The manuscript is organized as follows: Section 2 introduces the physical problem underlying active thin film modeling and highlights the main difficulties in analyzing such structures. Section 3 describes the numerical approach used to discretize the coupled electromechanical problem in space and time. Several numerical examples are proposed in Section 4 to demonstrate the capabilities of the presented approach. Specifically, we start simulating quasi-static and dynamic conditions to compare our approach with others proposed in the literature. Then we present some complex coupled electromechanical dynamic tests. Finally, Section 5 summarizes the results and highlights possible future studies.

## 2. Electromechanical model for thin-composite active films

MTFs (Shim et al., 2012; Böl et al., 2009) and similar constructs (MacQueen et al., 2018) are characterized, as a result of the manufacturing process, by two predominant spatial dimensions compared to the third one. Furthermore, the two materials naturally introduce a splitting of the body through the thickness: The complete domain $\Omega$ is composed of the active biological layer $\Omega_e$ of thickness $d_a$, and the elastic substrate $\Omega_s$ of thickness $d_s$ ($\Omega = \Omega_e \cup \Omega_s$ and $d = d_a + d_s$), as schematically represented in Fig. 1.

The electrophysiological activity in the cell layer is (substantially) uniform through the thickness of the biological tissue and unaffected by the substrate, while the mechanical response is more complex since the two materials, with completely different properties, interact. Indeed, a contraction of the cells generates an out-of-plane bending, resulting in a complex discontinuous state of stress, provided that the layers hold a perfect adhesion.

Numerical simulations of such a complex system are usually performed by considering a complete description of the construct, where every material layer defines a 3D subdomain (Shim et al., 2012). Conversely, we perform simpler simulations on two different bi-variate manifolds. We solve an electrophysiological excitation of the active layer, coupled to the mechanics of a Kirchhoff–Love shell, whose deformation is induced by the active layer contraction. These two sub-problems, numerically solved on their reference manifolds, are described in Sections 2.1 and 2.2, respectively.

### 2.1. Electrophysiological sub-problem

The electrophysiological sub-problem describes the cellular activity in response to external stimuli and the force generation responsible for tissue contraction at the local level. To represent such phenomena, different levels of schematizations can be used, ranging from uniform activations (Böl et al., 2009; Shim et al., 2012) to complex models based on differential equations (Göktepe and Kuhl, 2010; Pezzuto et al., 2014). The first type of approach applies in simplified situations, for instance, when the cell activity is uniform in space, while the latter can capture complex patterns of activation varying in space and time as a consequence of the interactions of multiple stimuli. In this work, we perform simulations using an activation field whose spatiotemporal distribution is described either via algebraic equations – the actual formulations are postponed to the examples in Section 4 – or by solving the so-called monodomain formulation, which is one of the classical models to simulate cardiac-like tissues (Franzone et al., 2014; Botti and Torre, 2023). This is herein described in its generic format, while the actual cell model and parameters used in the simulations are reported in Section 4.

The evolution in time $t \in [0, T]$ of the transmembrane potential $v$ in the subdomain $\Omega_e$ representing the active layer is given by the monodomain formulation. In our framework, every point $\mathbf{x}_e \in \Omega_e \subset \mathbb{R}^3$ is mapped using a reference 2D manifold $\bar{A}_e \subset \mathbb{R}^3$, that we identify as the mid-surface of the active layer highlighted in Fig. 1, and a third direction $\mathbf{h}_3$ locally orthogonal to the manifold:

$$\mathbf{x}_e = \mathbf{r}_e(\theta_e^1, \theta_e^2) + \theta_e^3 \, \mathbf{h}_3 \,, \tag{1}$$

where $\mathbf{r}_e(\theta_e^1, \theta_e^2)$ is a generic point of the manifold identified by the coordinates $\theta_e^1$, $\theta_e^2$ and $\theta_e^3 \in [-d_a/2, d_a/2]$ is the coordinate spanning the thickness of the layer $d_a$.

Assuming that the field $v$ is constant in the direction $\mathbf{h}_3$ because the active layer is thin, the following monodomain formulation can be solved on the reference manifold:

$$\begin{cases} C_m \dfrac{\partial v}{\partial t} = \nabla \cdot (\mathbf{D}\nabla v) - I^{ion} + I^{app} & \text{in } A_e \times [0, T] \\ \mathbf{n}_e \cdot \mathbf{D}\nabla v = I^n & \text{on } \partial A_e \times [0, T] \\ v(\mathbf{r}_e, 0) = v_0 & \text{in } \bar{A}_e \text{ for } t = 0 \end{cases} \tag{2}$$

being $C_m$ the capacitance of the cell membrane, $\mathbf{D}$ the tissue diffusivity tensor, $\mathbf{n}_e$ the outward-pointing vector on the boundary $\partial A_e$ of the surface, while $I^{app}$ and $I^n$ are source terms in the interior and on the boundary of the domain ($\bar{A}_e = A_e \cup \partial A_e$), respectively. Additionally, $v_0$ is the initial condition, usually representing the depolarized state of the cells.

We highlight that, in such equations, the operator $\nabla$ in (2) differentiates with respect to two directions defined in the local tangent plane.

Among the possible coordinate sets, we adopt a curvilinear framework[1] to describe the manifold, providing a clear identification of the local tangent plane for any generic manifold.

*Curvilinear reference frame.* Being $\mathbf{r}_e$ a point of the manifold defined by two coordinates $(\theta_e^1, \theta_e^2)$, the local tangent (covariant) vectors $\mathbf{h}_\alpha$ are defined (Itskov et al., 2007) as:

$$\mathbf{h}_\alpha = \frac{\partial \mathbf{r}_e}{\partial \theta_e^\alpha} = \mathbf{r}_{e,\alpha} \quad \text{for } \alpha = 1,2. \tag{3}$$

Accordingly, the normal vector is defined as:

$$\mathbf{h}_3 = \frac{\mathbf{h}_1 \times \mathbf{h}_2}{\|\mathbf{h}_1 \times \mathbf{h}_2\|}, \tag{4}$$

completing the definition of the curvilinear reference frame. In this work, Latin indices take on values 1, 2, 3 (if not specified differently), while Greek indices take on values 1, 2. Einstein's indicial notation is used (repeated indices imply summation unless differently specified).

The potential on the manifold is a function of the coordinates $v = v(\theta_e^\alpha)$, consequently the non-null components of its gradient lie on the local tangent plane and read:

$$\nabla v = v_{,\alpha} \mathbf{h}^\alpha, \tag{5}$$

being the contravariant vectors $\mathbf{h}^\alpha$ computed in terms of covariant vectors and of the inverse of the metric tensor $h_{\alpha\beta}$ as:

$$\mathbf{h}^\alpha = h^{\alpha\beta} \mathbf{h}_\beta, \tag{6}$$

where

$$h^{\alpha\beta} = \left[h_{\alpha\beta}\right]^{-1} = \left[\mathbf{h}_\alpha \cdot \mathbf{h}_\beta\right]^{-1}. \tag{7}$$

Note that, according to our assumptions, $\mathbf{h}^3 = \mathbf{h}_3$.

*Mondomain formulation in curvilinear coordinates.* Assuming that the tissue conductivity is constant and isotropic ($\mathbf{D} = D\mathbf{I}$), (2) can be expressed in curvilinear coordinates as:

$$\begin{cases} C_m \dfrac{\partial v}{\partial t} = D\Delta v - I^{ion} + I^{app} & \text{in } A_e \times [0,T] \\ \mathbf{n}_e \cdot D\, v_{,\alpha} \mathbf{h}^\alpha = I^n & \text{on } \partial A_e \times [0,T] \\ v(\mathbf{r}_e, 0) = v_0 & \text{in } \bar{A}_e \text{ for } t = 0 \end{cases} \tag{8}$$

where

$$\Delta v = \left(v_{,\alpha\beta} \mathbf{h}^\alpha - v_{,\alpha} \Gamma_{\beta k}^\alpha \mathbf{h}^k\right) \cdot \mathbf{h}^\beta, \tag{9}$$

and $\Gamma_{\beta k}^\alpha$ is the Christoffel symbol of second kind (Itskov et al., 2007).

*Cell model.* The term $I^{ion}$ in the previous equation represents the ionic current through the cell membrane and is a nonlinear term whose amplitude depends on a set of state variables $\mathbf{w}$ – which represent the concentration of ions and the activation state of ionic channels in the membrane – evolving in time:

$$\begin{cases} I^{ion} = f(v, \mathbf{w}) & \text{in } \bar{A}_e \times [0,T] \\ \dfrac{\partial \mathbf{w}}{\partial t} = \mathbf{f}_w(v, \mathbf{w}) & \text{in } \bar{A}_e \times [0,T] \end{cases}, \tag{10}$$

We note that, since the cell model is point-wise defined (i.e., the ionic current and the state variable in a given point depend on quantities defined at that specific point only) and no spatial derivatives are involved, such equations read the same in cartesian and curvilinear coordinates.

*Force generation.* Finally, the evolution of the active stress $\sigma_a$ generated by a cell depends on the potential $v$, possibly on $\mathbf{w}$, and is given by an additional set of differential equations as follows:

$$\begin{cases} \dfrac{\partial \sigma_a}{\partial t} = q\left(\sigma_a, v, \mathbf{w}, \mathbf{s}\right) & \text{in } \bar{A}_e \times [0,T] \\ \dfrac{\partial \mathbf{s}}{\partial t} = \mathbf{q}_s\left(\sigma_a, v, \mathbf{w}, \mathbf{s}\right) & \text{in } \bar{A}_e \times [0,T] \end{cases}, \tag{11}$$

where $\mathbf{s}$ is a second set of state variables representing the cell activity.

In our notation, $f$, $\mathbf{f}_w$, $q$, and $\mathbf{q}_s$ are scalar and vector functions characterizing the specific type of cell under investigation. For the actual models herein simulated and the initial conditions, the reader is referred to Section 4.

In the present work, for the sake of simplicity, we solve the mondomain formulation assuming that the domain is rigid, neglecting mechano-electrical feedback and the effect of the tissue deformations on the diffusivity tensor. However, they can be included in the formulation as shown in previous works (Nitti et al., 2021).

### 2.2. Mechanical sub-problem

In the mechanical sub-problem, differently from the electrophysiological one, both the active layer and the elastic substrate are explicitly modeled since the hypothesis of perfect adhesion (Shim et al., 2012; Pezzuto et al., 2014) couples the domains.

In similarity to (1), a material point $\dot{\mathbf{x}} \in \Omega$ is mapped as:

$$\dot{\mathbf{x}} = \dot{\mathbf{r}}\left(\theta^1, \theta^2\right) + \theta^3 \hat{\mathbf{a}}_3, \tag{12}$$

where $\dot{\mathbf{r}}$ is a point of the shell mid-surface $\bar{A} = A \cup \partial A$, $\hat{\mathbf{a}}_3$ is the normal to the manifold, and $\theta^3 \in [-d/2, d/2]$ the corresponding coordinate in the normal direction. Since the domain deforms under the effect of external actions, we use the symbol $\dot{\bullet}$ to identify a generic quantity referred to the reference configuration.

In the reference configuration, a point in the active layer can be mapped by combining (1) and (12):

$$\dot{\mathbf{x}}_e = \dot{\mathbf{r}} + \frac{(d - d_a)}{2} \hat{\mathbf{a}}_3 + \theta_e^3 \mathbf{h}_3. \tag{13}$$

Moreover, we assume that the two normal vectors are co-axial (i.e., $\mathbf{h}_3 = \hat{\mathbf{a}}_3$). Therefore, (13) reduces to:

$$\dot{\mathbf{x}}_e = \dot{\mathbf{r}}\left(\theta^1, \theta^2\right) + \left(\theta_e^3 + \frac{d - d_a}{2}\right) \hat{\mathbf{a}}_3 \tag{14}$$

for $\theta_e^3 + (d - d_a)/2 = \theta^3 \in \left[d/2 - d_a, d/2\right]$, therefore the upper surface of the active layer $\theta^3 = d/2$ corresponds to the coordinate $\theta_e^3 = d_a/2$, for instance. In the limits of $d_a \to d$ or $d_a \to 0$, we recover the case of a shell constituted solely by active or elastic material, respectively.

In mechanics, the Kirchhoff–Love shell theory (Kiendl et al., 2015) is a well-established approximation for thin structures. In this work, we rely on such an approach to model MTFs, following the developments of previous works (Nitti et al., 2021; Kiendl et al., 2015), herein summarized for the sake of completeness.

*Kinematics.* In the deformed configuration, a material point is related to its reference configuration by the following relation:

$$\mathbf{x} = \mathbf{r} + \theta^3 \mathbf{a}_3 = \dot{\mathbf{r}} + \mathbf{u} + \theta^3 \mathbf{a}_3, \tag{15}$$

being $\mathbf{u}(\theta^1, \theta^2)$ the displacement field on the shell mid-surface, which represents the unknown of the mechanical problem.

As done in the electrophysiological problem, we define a curvilinear reference frame on the mid-plane of the shell:

$$\mathbf{a}_\alpha = \mathbf{r}_{,\alpha} = \dot{\mathbf{r}}_{,\alpha} + \mathbf{u}_{,\alpha}, \tag{16}$$

while the vector normal to the surface reads:

$$\mathbf{a}_3 = \frac{\mathbf{a}_1 \times \mathbf{a}_2}{\|\mathbf{a}_1 \times \mathbf{a}_2\|}. \tag{17}$$

---

[1] The interested reader is referred to Appendix A and Itskov et al. (2007) for more details on curvilinear reference frames.

Analogously, deriving (15), the covariant vectors for points lying out of the mid-surface (i.e., $\theta^3 \neq 0$) assume the following forms:

$$\mathbf{g}_\alpha = \mathbf{a}_\alpha + \theta^3 \mathbf{a}_{3,\alpha} \,, \tag{18}$$

while

$$\mathbf{g}_3 = \mathbf{a}_3 \,. \tag{19}$$

Similar relations hold for the reference configuration as well, and are used in the following derivations to characterize the stresses and strains in the domain.

In the finite elasticity setting, the Cauchy–Green deformation tensor

$$\mathbf{C} = C_{ij} \, \mathring{\mathbf{g}}^i \otimes \mathring{\mathbf{g}}^j \,, \tag{20}$$

whose components – according to the Kirchhoff–Love kinematics (Kiendl et al., 2015) – are:

$$C_{ij} = \begin{pmatrix} g_{11} & g_{12} & 0 \\ g_{21} & g_{22} & 0 \\ 0 & 0 & C_{33} \end{pmatrix} \,, \tag{21}$$

is indeed a function of the metric tensor

$$g_{\alpha\beta} = \mathbf{g}_\alpha \cdot \mathbf{g}_\beta \,, \tag{22}$$

and of the contravariant vectors $\mathring{\mathbf{g}}^i$. Note that, the through-the-thickness component $C_{33}$ depends on the in-plane strains and is computed such that plane-stress constraints are fulfilled (Kiendl et al., 2015). In shells, only the in-plane components $C_{\alpha\beta}$ actively contribute to the internal virtual work via the Green–Lagrange strain tensor $\mathbf{E}$:

$$E_{\alpha\beta} = \frac{1}{2} \left( C_{\alpha\beta} - \mathring{g}_{\alpha\beta} \right) \,. \tag{23}$$

which is work-conjugate with the second Piola–Kirchhoff stress tensor.

To relate the in-plane components of the Green–Lagrange tensor to the displacements of the mid-surface, we introduce the membrane strains $\varepsilon_{\alpha\beta}$ and the curvatures $\kappa_{\alpha\beta}$ as:

$$\varepsilon_{\alpha\beta} = \frac{1}{2} \left( a_{\alpha\beta} - \mathring{a}_{\alpha\beta} \right) \tag{24}$$

and

$$\kappa_{\alpha\beta} = \mathring{b}_{\alpha\beta} - b_{\alpha\beta} \,, \tag{25}$$

where $\mathbf{u}$ enters through $\mathbf{a}_\alpha$ ($\mathbf{a}_\beta$) in the following definitions:

$$a_{\alpha\beta} = \mathbf{a}_\alpha \cdot \mathbf{a}_\beta \,, \quad b_{\alpha\beta} = \mathbf{a}_{\alpha,\beta} \cdot \mathbf{a}_3 \,. \tag{26}$$

Accordingly, the Green–Lagrange strain components are expressed as:

$$E_{\alpha\beta} = \varepsilon_{\alpha\beta} + \theta^3 \kappa_{\alpha\beta} \,, \tag{27}$$

and can be used to compute the stress tensor and the internal energy.

*Principle of virtual work.* According to the principle of virtual works, the solution of the mechanical problem $\mathbf{u}$ is such that the total energy $\mathcal{W}$ is minimized:

$$\delta \mathcal{W} = \delta \mathcal{W}_\rho + \delta \mathcal{W}_{int} - \delta \mathcal{W}_{ext} = 0 \quad \forall \delta \mathbf{u} \,, \tag{28}$$

being $\delta \mathcal{W}$ the energy variation with respect to any admissible virtual displacement field $\delta \mathbf{u}$. In a Kirchhoff–Love shell, the total energy variation is given by the virtual work of inertia forces $\delta \mathcal{W}_\rho$ and internal actions $\delta \mathcal{W}_{int}$ as well as by the of virtual work of the external loads $\delta \mathcal{W}_{ext}$. The first two terms read:

$$\delta \mathcal{W}_\rho = \int_A \rho \ddot{\mathbf{u}} \cdot \delta \mathbf{u} \, dA \,, \tag{29}$$

and

$$\delta \mathcal{W}_{int} = \int_A \mathbf{n} : \delta \boldsymbol{\varepsilon} + \mathbf{m} : \delta \boldsymbol{\kappa} \, dA \,. \tag{30}$$

Any virtual displacement induces a variation in the membrane strain $\delta \boldsymbol{\varepsilon}$ and curvatures $\delta \boldsymbol{\kappa}$, which are computed by taking the Gateaux derivative with respect to $\delta \mathbf{u}$ of (24) and (25) and are work-conjugated to the normal actions and bending moments $\mathbf{n}$ and $\mathbf{m}$, respectively. Moreover, in (29), the product of the density per unit of area $\rho$ and acceleration $\ddot{\mathbf{u}} = \frac{\partial^2 \mathbf{u}}{\partial t^2}$ represents the inertia forces, which are neglected in quasi-static simulations. In the dynamic case, instead, initial displacements $\mathbf{u}$ and velocities $\dot{\mathbf{u}}$ are prescribed at time $t = 0$ for the mid-surface to complement the mathematical formulation.

The external work $\delta \mathcal{W}_{ext}$ is given by:

$$\delta \mathcal{W}_{ext} = \int_A \mathbf{p} \cdot \delta \mathbf{u} \, dA + \int_{\partial A_N} \mathbf{t} \cdot \delta \mathbf{u} \, ds \,. \tag{31}$$

and in the numerical test presented in this work is assumed to be null since the body is unloaded.

We recall that, to be admissible, a virtual displacement must be null on the boundary of the manifold $\partial A_D$ where Dirichlet boundary conditions $\mathbf{u} = \mathbf{u}_D$ are specified. On the remaining part of the boundary $\partial A_N$, tractions $\mathbf{t}$ are prescribed ($\partial A = \partial A_N \cup \partial A_D$ and $\partial A_N \cap \partial A_D = \emptyset$), while body loads $\mathbf{p}$ are prescribed in the interior of the manifold. Moreover, in computing such integrals, $dA$ represents the differential area in the reference configuration:

$$dA = \sqrt{\det(\mathring{a}_{\alpha\beta})} \, d\theta^1 d\theta^2 \,, \tag{32}$$

while the differential length in the direction spanned by the coordinate $\theta^\alpha$ is

$$ds = \|\mathring{\mathbf{g}}_\alpha\| \, d\theta^\alpha \,. \tag{33}$$

Finally, we define the normal force $\mathbf{n}$ and bending moment $\mathbf{m}$ in (30) via integration through-the-thickness of the second Piola–Kirchhoff stress tensor:

$$n^{\alpha\beta} = \int_{-d/2}^{d/2} S^{\alpha\beta} \, d\theta^3 \,, \tag{34}$$

$$m^{\alpha\beta} = \int_{-d/2}^{d/2} S^{\alpha\beta} \theta^3 \, d\theta^3 \,. \tag{35}$$

Differently from previous works on shells, the stress tensor in the active layer is assumed to directly depend not only on the strains but on the electrophysiological activation as well.

*Active stress formulation.* The electrophysiological activation is modeled according to the active stress approach: The second Piola–Kirchhoff stress tensor is given by the sum of the classical elastic term, given by the derivative of the strain energy function $\psi$, and by the so-called active term $\mathbf{S}_a$ (Ambrosi and Pezzuto, 2012; Sundnes et al., 2014), as follows:

$$\mathbf{S} = \frac{\partial \psi}{\partial \mathbf{E}} + \mathbf{S}_a = 2 \frac{\partial \psi}{\partial \mathbf{C}} + \mathbf{S}_a \,. \tag{36}$$

In this work, we assume that contractile cells are organized in fibers oriented in the plane (i.e., out-of-plane stress components are null) and that $\mathbf{S}_a$ is independent of the transversal strain, although generalizations are possible. Moreover, we assume that the passive response of the tissue is stiffer in the fiber direction enforcing an anisotropic formulation for the strain energy $\psi$ (i.e., the energy is a function of the fiber orientation).

Given the stress amplitude $\sigma_a$ predicted by the activation model, the stress tensor, after a pullback on the reference configuration, reads:

$$\mathbf{S}_a = \frac{\sigma_a}{\lambda_f^2} \left( \mathring{\mathbf{f}}_0 \otimes \mathring{\mathbf{f}}_0 \right) \,, \tag{37}$$

where $\mathring{\mathbf{f}}_0$ is the vector representing the fiber direction, expressed in cartesian coordinates, and $\lambda_f$ is the stretch in the fiber direction. In curvilinear coordinates, the in-plane components read:

$$\mathbf{S}_a = S_a^{\alpha\beta} \mathring{\mathbf{g}}_\alpha \otimes \mathring{\mathbf{g}}_\beta \quad \text{with} \quad S_a^{\alpha\beta} = \frac{\sigma_a}{\lambda_f^2} \, \mathring{f}_0^\alpha \mathring{f}_0^\beta \,, \tag{38}$$

being

$$\mathring{\mathbf{f}}_0^\alpha = \mathring{\mathbf{f}}_0 \cdot \mathring{\mathbf{g}}^\alpha \quad \text{and} \quad \lambda_f = \sqrt{\mathring{\mathbf{f}}_0^\alpha \, C_{\alpha\beta} \, \mathring{\mathbf{f}}_0^\beta} \,. \tag{39}$$

In the elastic substrate, the active component of the stress is null. We highlight that the passive response of the two materials constituting the active layer and substrate are completely different. For instance, the shear modulus of the substrate is in the order of a few hundred kPa, while for the biological tissue it is in the order of one kPa (Shim et al., 2012). This is reflected in a discontinuity in the in-plane stress distribution through the thickness (see Fig. 1).

## 3. Discrete formulation

In presenting the discrete version of the model introduced in Section 2, we revert the presentation strategy. The mechanical problem is discretized first in Section 3.1, then the dynamics is introduced in Section 3.2, such that the spatial and temporal discretizations of the electrophysiological problem, in Sections 3.3 and 3.4, and the coupling scheme, in Section 3.5, are then built on that basis. The same structure is adopted in presenting the numerical examples to better focus on the various features of the algorithms.

In Kirchhoff–Love shells, the curvature depends on high-order derivatives, as shown in (26), which implies requiring at least $C^1$-continuity across the elements for the basis functions to enable a discretization in primal form. Quadratic (or higher order) B-spline discretizations easily satisfy such a requirement. Furthermore, B-splines naturally adapt to the curvilinear framework thanks to the concept of parametric space: Given a parametric square spanned by the coordinates $\theta^\alpha$, where the B-splines are defined, every point is mapped into the shell mid-surface via a linear combination of control point coordinates $\mathbf{B}$ and basis functions as follows:

$$\mathring{r}_j = N_{ji}(\theta^1, \theta^2)\mathbf{B}_i \quad \text{for } i = 1, \dots, 3 \times n \,, \tag{40}$$

where the $n$ B-spline functions are repeated in the matrix $\mathbf{N}$ to represent a vector field, as done in standard implementations of the Galerkin method. The non-identically-null functions are defined via tensor products of univariate B-splines using the Cox-De Boor formula, given the polynomial degree $p_a$ and the open knot vector $\Xi^\alpha = \{\xi_1^\alpha, \xi_2^\alpha, \dots, \xi_{m_a + p_a + 1}^\alpha\}$ – where $\xi_i^\alpha$ is the $i$th knot – per every parametric direction (i.e., $\alpha = 1$ or $\alpha = 2$ and $n = m_1 \times m_2$). Moreover, a similar definition holds for the reference surface of the active layer, where a new set of control point coordinates $\mathbf{B}_c$ is computed shifting the shell mid-plane in the normal direction. The mechanical and electrophysiological sub-problems are solved using such splines, possibly refined via knot-insertion and degree-elevation algorithms to achieve a suitable level of accuracy in each sub-problem independently.

### 3.1. Spatial discretization for active Kirchhoff–Love shells

Exploiting the isoparametric concept, the displacement field $\mathbf{u}$ is approximated as:

$$\mathbf{u} = \mathbf{N}(\theta^1, \theta^2)\hat{\mathbf{u}}(t) \,, \tag{41}$$

being $\hat{\mathbf{u}}$ a column-vector collecting the control displacements, whose values vary in time. Accordingly, the deformed mid-surface position results in

$$\mathbf{r} = \mathbf{N}(\mathbf{B} + \hat{\mathbf{u}}) \tag{42}$$

and the discrete definition of the covariant vectors and derivatives read:

$$\mathbf{a}_\alpha = \mathbf{N}_{,\alpha}(\mathbf{B} + \hat{\mathbf{u}}) \,, \quad \mathbf{a}_{\alpha,\beta} = \mathbf{N}_{,\alpha\beta}(\mathbf{B} + \hat{\mathbf{u}}) \,. \tag{43}$$

An expansion of the virtual displacements similar to (41)

$$\delta\mathbf{u} = \mathbf{N}\,\delta\hat{\mathbf{u}} \tag{44}$$

leads to the following residual equation:

$$\mathbf{R}(\hat{\mathbf{u}}, \ddot{\hat{\mathbf{u}}}) = \mathbf{M}\ddot{\hat{\mathbf{u}}} + \mathbf{F}_{int} - \mathbf{F}_{ext} = \mathbf{0} \quad \forall \delta\hat{\mathbf{u}} \,, \tag{45}$$

where

$$\mathbf{M} = \int_A \mathbf{N}^T \rho \mathbf{N} \, dA \,, \tag{46}$$

and $\mathbf{F}_{int}$ and $\mathbf{F}_{ext}$ are the discrete version of (30) and (31), respectively. In the present work, (45) is solved at certain time instants, specified in Section 3.2, via Newton method, which requires the linearization $\mathbf{K}$ of the residual with respect to the control variables $\hat{\mathbf{u}}$:

$$\mathbf{K} = D_{\hat{\mathbf{u}}}[\mathbf{R}] = \mathbf{K}_{int} - \mathbf{K}_{ext} \,, \tag{47}$$

being

$$\mathbf{K}_{int} = \int_A D_{\hat{\mathbf{u}}}[\mathbf{n}] : \delta\boldsymbol{\varepsilon} + \mathbf{n} : D_{\hat{\mathbf{u}}}[\delta\boldsymbol{\varepsilon}] \, dA$$
$$+ \int_A D_{\hat{\mathbf{u}}}[\mathbf{m}] : \delta\boldsymbol{\kappa} + \mathbf{m} : D_{\hat{\mathbf{u}}}[\delta\boldsymbol{\kappa}] \, dA \tag{48}$$

and $\mathbf{K}_{ext}$ the linearization of the external loads, that in the case of non-follower external actions is null.

In this work, the active stress modifies the linearization of the internal actions compared to classical shell formulations. In fact, the linearization of the normal forces:

$$D_{\hat{\mathbf{u}}}[n^{\alpha\beta}] = \int_{-d/2}^{d/2} \hat{\mathbb{C}}^{\alpha\beta\gamma\delta} d\theta^3 \, D_{\hat{\mathbf{u}}}[\varepsilon_{\gamma\delta}]$$
$$+ \int_{-d/2}^{d/2} \hat{\mathbb{C}}^{\alpha\beta\gamma\delta} \theta^3 d\theta^3 \, D_{\hat{\mathbf{u}}}[\kappa_{\gamma\delta}] \tag{49}$$

and of the bending moments:

$$D_{\hat{\mathbf{u}}}[m^{\alpha\beta}] = \int_{-d/2}^{d/2} \hat{\mathbb{C}}^{\alpha\beta\gamma\delta} \theta^3 d\theta^3 \, D_{\hat{\mathbf{u}}}[\varepsilon_{\gamma\delta}]$$
$$+ \int_{-d/2}^{d/2} \hat{\mathbb{C}}^{\alpha\beta\gamma\delta} (\theta^3)^2 d\theta^3 \, D_{\hat{\mathbf{u}}}[\kappa_{\gamma\delta}] \tag{50}$$

depend on $\hat{\mathbb{C}}^{\alpha\beta\gamma\delta}$, that is the linearization of $\mathbf{S}$ with respect to the strains in plane-stress conditions, which is given by the standard contribution of the term derived by the strain energy $\hat{\mathbb{C}}_p^{\alpha\beta\gamma\delta}$ (see Kiendl et al., 2015) plus the active term:

$$\hat{\mathbb{C}}^{\alpha\beta\gamma\delta} = \hat{\mathbb{C}}_p^{\alpha\beta\gamma\delta} + 2 \frac{\partial S_a^{\alpha\beta}}{\partial C_{\gamma\delta}} \,. \tag{51}$$

Note that, in this last equation, the hypothesis on the active stress acting only in the plane is enforced and that integrals (49), (50) and the discrete version of (34) and (35) are computed numerically distributing uni-variate Gauss quadrature points through the thickness of every material layer. When a different number of layers is considered, the approach simply requires the subdivision of the thickness into more layers during integration, without any further modification.

In the case of quasi-static problems, the inertia terms are neglected in the equilibrium equation and a generic iteration of the Newton method reads:

$$\mathbf{K}\Delta\hat{\mathbf{u}} = -\mathbf{R} \,, \tag{52}$$

being $\Delta\hat{\mathbf{u}}$ the increment to be added to the control displacements. Note that an arc-length method (Verhelst et al., 2024) may be required to solve this problem, for instance, if structural buckling occurs.

Verifications of the proposed approach are reported in Appendix B.

### 3.2. Temporal discretization for active-shell dynamics

When dynamics is considered, the time is discretized in intervals, herein assumed of equal size $\Delta t$, such that $t_{i+1} = t_i + \Delta t$ for $i = 1, \dots, T/\Delta t$, and the residual (45) is solved within each interval $[t_i, t_{i+1}]$ following the generalized-$\alpha$ method.

At $\alpha$-time, the displacement, velocity, and acceleration control variables can be expressed in terms of values at time $t_i$ and $t_{i+1}$ as:

$$\hat{\mathbf{u}}_\alpha = \alpha_f \hat{\mathbf{u}}_{i+1} + (1 - \alpha_f) \hat{\mathbf{u}}_i, \tag{53}$$

$$\dot{\hat{\mathbf{u}}}_\alpha = \alpha_f \dot{\hat{\mathbf{u}}}_{i+1} + (1 - \alpha_f) \dot{\hat{\mathbf{u}}}_i, \tag{54}$$

$$\ddot{\hat{\mathbf{u}}}_\alpha = \alpha_m \ddot{\hat{\mathbf{u}}}_{i+1} + (1 - \alpha_m) \ddot{\hat{\mathbf{u}}}_i, \tag{55}$$

while the updates for displacements and velocities are given by the Newmark formulas:

$$\hat{\mathbf{u}}_{i+1} = \hat{\mathbf{u}}_i + \Delta t \dot{\hat{\mathbf{u}}}_i + \frac{\Delta t^2}{2} \left( (1 - 2\beta) \ddot{\hat{\mathbf{u}}}_i + 2\beta \ddot{\hat{\mathbf{u}}}_{i+1} \right) \tag{56}$$

$$\dot{\hat{\mathbf{u}}}_{i+1} = \dot{\hat{\mathbf{u}}}_i + \Delta t \left( (1 - \gamma) \ddot{\hat{\mathbf{u}}}_i + \gamma \ddot{\hat{\mathbf{u}}}_{i+1} \right). \tag{57}$$

The residual equation:

$$\mathbf{R}_\alpha = \mathbf{M} \ddot{\hat{\mathbf{u}}}_\alpha + \mathbf{F}_{int} \left( \hat{\mathbf{u}}_\alpha \right) - \mathbf{F}_{ext} = \mathbf{0}, \tag{58}$$

can be solved via the Newton method linearizing with respect to the accelerations:

$$\frac{d\mathbf{R}_\alpha}{d\ddot{\hat{\mathbf{u}}}_{i+1}} \Delta \ddot{\mathbf{u}}_{i+1} = -\mathbf{R}_\alpha, \tag{59}$$

which leads to the following linear system:

$$\frac{d\mathbf{R}_\alpha}{d\ddot{\hat{\mathbf{u}}}_{i+1}} \Delta \ddot{\mathbf{u}}_{i+1} = -\mathbf{M} \ddot{\mathbf{u}}_\alpha - \mathbf{F}_\alpha^{int} + \mathbf{F}_\alpha^{ext}, \tag{60}$$

being

$$\frac{d\mathbf{R}_\alpha}{d\ddot{\hat{\mathbf{u}}}_{i+1}} = \alpha_m \mathbf{M} + \alpha_f \beta \Delta t^2 \mathbf{K}_\alpha. \tag{61}$$

In the specific case of damping proportional to the velocity, as in the case of Rayleigh model for viscous damping, the formulation is modified accordingly as:

$$\mathbf{R}_\alpha = \mathbf{M} \ddot{\hat{\mathbf{u}}}_\alpha + \mathbf{C} \dot{\hat{\mathbf{u}}}_\alpha + \mathbf{F}_{int} \left( \hat{\mathbf{u}}_\alpha \right) - \mathbf{F}_{ext} = \mathbf{0}, \tag{62}$$

being $\mathbf{C}$ the damping matrix, and the resulting linear system is:

$$\frac{d\mathbf{R}_\alpha}{d\ddot{\hat{\mathbf{u}}}_{i+1}} \Delta \ddot{\mathbf{u}}_{i+1} = -\mathbf{M} \ddot{\mathbf{u}}_\alpha - \mathbf{C} \dot{\hat{\mathbf{u}}}_\alpha - \mathbf{F}_\alpha^{int} + \mathbf{F}_\alpha^{ext}, \tag{63}$$

where the left-hand side matrix reads as:

$$\frac{d\mathbf{R}_\alpha}{d\ddot{\hat{\mathbf{u}}}_{i+1}} = \alpha_m \mathbf{M} + \alpha_f \gamma \Delta t \mathbf{C} + \alpha_f \beta \Delta t^2 \mathbf{K}_\alpha. \tag{64}$$

In the present work, we adopt classical values for the parameters characterizing the time integration scheme:

$$\alpha_m = \frac{2 - \rho_\infty}{1 + \rho_\infty}, \ \alpha_f = \frac{1}{1 + \rho_\infty}, \tag{65}$$

$$\beta = \frac{(1 - \alpha_f + \alpha_m)^2}{4}, \ \gamma = \frac{1}{2} - \alpha_f + \alpha_m, \tag{66}$$

and

$$\rho_\infty = 0.5. \tag{67}$$

Finally, the internal work $\mathbf{F}_\alpha^{int}$ depends on the level of mechanical activation of cells at time $t_\alpha$. We assume that the time dependence of such a field is given by linear interpolation of the values at the extremities of the time interval:

$$\sigma_\alpha(\hat{\mathbf{x}}, t_\alpha) = \alpha_f \sigma_\alpha(\hat{\mathbf{x}}, t_i) + (1 - \alpha_f) \sigma_\alpha(\hat{\mathbf{x}}, t_{i+1}), \tag{68}$$

in similarity to (53).

### 3.3. Spatial discretization of the electrophysiological sub-problem

The shape of the shell and the related mechanical problem are discretized via a spline-based approach and we adopt the same type of splines to solve the electrophysiological problem. However, instead of using an isogeometric-Galerkin approach, we adopt an isogeometric-Collocation method, as proposed in Torre et al. (2022, 2023a). Part of the computational effort in solving electrophysiology is entailed by the discretization of the reactive term ($I^{ion}$), which is coupled to the system of equations describing the cell activity (10). The collocation approach aims at limiting the computational effort avoiding the quadrature of such a term in favor of the evaluation of the reactive term in a minimal set of points. It requires basis functions that are at least $C^1$-continuous, but this is anyway the case for the splines used in the mechanical sub-problem.

According to the collocation method, the monodomain formulation (8) is discretized by evaluating the strong form of the differential problem at every $j$th collocation point $\bar{\mathbf{r}}_{ej}$, defined as the image of the corresponding Greville abscissa via a standard linear combination of B-splines and control point coordinates of the reference configuration:

$$\bar{\mathbf{r}}_{ej} = \mathbf{N}^{ge} \left( \bar{\boldsymbol{\xi}}_j \right) \mathbf{B}_e. \tag{69}$$

Following this approach, the bi-variate coordinates $\bar{\boldsymbol{\xi}}_j = (\bar{\xi}_i^1, \bar{\xi}_k^2)$ are defined via the tensor product of uni-variate abscissae:

$$\bar{\xi}_i^\alpha = \frac{\xi_{i+1}^\alpha + \cdots + \xi_{i+p_\alpha}^\alpha}{p_\alpha} \quad \forall i = 1, \ldots, m_\alpha, \tag{70}$$

and, therefore, are in the number of the control points (i.e., $j = 1, \ldots, n_e$).

Introducing an approximation similar to (41) for the scalar potential:

$$\upsilon = \mathbf{N}^e \left( \theta^1, \theta^2 \right) \hat{\mathbf{v}}(t), \tag{71}$$

where $\mathbf{N}^e$ is the row-vector collecting only the non-identically-null basis functions of the matrix $\mathbf{N}^{ge}$, (8) is discretized by evaluations at the collocation points in the interior of the domain $A_e$:

$$\mathbf{M}_j^{ei} \dot{\hat{\mathbf{v}}} = \mathbf{K}_j^{ei} \hat{\mathbf{v}} - \mathbf{I}_j^i + \mathbf{I}_j^{app} \quad \forall j = 1, \ldots, n_i, \tag{72}$$

being $n_i = n + 4 - (m_1 \times 2) - (m_2 \times 2)$ the number of control points in the interior of the domain:

$$\mathbf{I}_j^i = I^{ion} \left( \upsilon(\bar{\mathbf{r}}_{ej}), t \right), \quad \mathbf{I}_j^{app} = I^{app} \left( \bar{\mathbf{r}}_{ej}, t \right), \tag{73}$$

and

$$\mathbf{M}_j^{ei} = C_m \mathbf{N}^e \left( \bar{\boldsymbol{\xi}}_j \right), \tag{74}$$

$$\mathbf{K}_j^{ei} = D \mathbf{h}^\beta \left( \mathbf{h}^\alpha \mathbf{N}_{,\alpha\beta}^e \left( \bar{\boldsymbol{\xi}}_j \right) - \Gamma_{\beta k}^\alpha \mathbf{h}^k \mathbf{N}_{,\alpha}^e \left( \bar{\boldsymbol{\xi}}_j \right) \right). \tag{75}$$

Herein, the Christoffel symbol and curvilinear frames are evaluated by exploiting the geometrical mapping, which is known.

Similarly, Neumann boundary conditions are evaluated at collocation points on the boundary of the parametric space defining $n_b = (m_1 \times 2) + (m_2 \times 2) - 4$ equations, for which the $j$th reads:

$$\mathbf{L}_j^{eb} \mathbf{v} = \mathbf{I}_j^{nb}. \tag{76}$$

In such an equation, the left-hand-side is defined as:

$$\mathbf{L}_j^{eb} = D \mathbf{a}_e^\alpha \cdot \mathbf{n}_e \mathbf{N}_{,\alpha}^e \left( \bar{\boldsymbol{\xi}}_j \right), \tag{77}$$

and the right-hand side as:

$$\mathbf{I}_j^{nb} = I^n \left( \bar{\mathbf{r}}_{ej}, t \right). \tag{78}$$

As a remark, we note that the derivation herein presented corresponds to a Petrov–Galerkin formulation in which the test functions are Dirac Delta functions.

Finally, the $n_e = n_i + n_b$ semi-discrete equations can be collected in matrix form for notation convenience as follows:

$$(\mathbf{M}^e + \mathbf{L}^e)\,\dot{\hat{\mathbf{v}}} = \mathbf{K}^e\,\hat{\mathbf{v}} - \mathbf{I}^{ion} + \mathbf{I}^{app} + \mathbf{I}^n\,, \tag{79}$$

where the first row corresponds to the discretization of the domain interior and the second one to the boundary:

$$\mathbf{M}^e = \begin{bmatrix} \mathbf{M}^{ei} \\ \mathbf{0} \end{bmatrix}\,, \quad \mathbf{L}^e = \begin{bmatrix} \mathbf{0} \\ \mathbf{L}^{eb} \end{bmatrix}\,, \quad \mathbf{K}^e = \begin{bmatrix} \mathbf{K}^{ei} \\ \mathbf{0} \end{bmatrix}\,, \tag{80}$$

and

$$\mathbf{I}^{ion} = \begin{bmatrix} \mathbf{I}^i \\ \mathbf{0} \end{bmatrix}\,, \quad \mathbf{I}^{app} = \begin{bmatrix} \mathbf{I}^{app} \\ \mathbf{0} \end{bmatrix}\,, \quad \mathbf{I}^n = \begin{bmatrix} \mathbf{0} \\ \mathbf{I}^{nb} \end{bmatrix}\,. \tag{81}$$

Verifications of the proposed approach are reported in Appendix C.

### 3.4. Time discretization for the electrophysiological sub-problem

To discretize (79) in time, we adopt a standard semi-implicit approach based on a combination of the Crank–Nicholson scheme and the two-step Adam–Bashforth method, as previously done in Nitti et al. (2023). Subdividing the time into equally spaced intervals of amplitude $\Delta t_e = \Delta t/k$ (being $k$ a positive integer) a generic integration step to compute the potential at the subsequent time point, $\mathbf{v}_{i+1}$, reads:

$$\mathbf{A}\hat{\mathbf{v}}_{i+1} = \mathbf{b}^e_{i+1}\,, \tag{82}$$

being

$$\mathbf{A} = \mathbf{M}^e - \frac{\Delta t_e}{2}\mathbf{K}^e + \mathbf{L}^e\,, \tag{83}$$

and

$$\begin{aligned} \mathbf{b}^e_{i+1} = \tilde{\mathbf{v}}_i + \frac{\Delta t_e}{2}\mathbf{K}^e\hat{\mathbf{v}}_i + \frac{\Delta t_e}{2}\left(\mathbf{I}^{app}_{i+1} + \mathbf{I}^{app}_i\right) \\ + \frac{3\Delta t_e}{2}\mathbf{I}^{ion}_i - \frac{\Delta t_e}{2}\mathbf{I}^{ion}_{i-1} + \mathbf{I}^n_{i+1}\,. \end{aligned} \tag{84}$$

We note that, for the first time step, a single-step Adam–Bashforth scheme is used.

In (84), $\tilde{\mathbf{v}}_i$ is given by the product of $\mathbf{M}^e$ and control variables $\hat{\mathbf{v}}_i$ and it represents, up to the multiplicative constant $C_m$, the transmembrane potential at collocation points in the internal part of the domain $A_e$. Such values are used to integrate in time the cell model (10), via an explicit fourth-order Runge Kutta method, and to compute the ionic currents at these specific locations.

### 3.5. One-way coupling combining Galerkin- and collocation-based solvers

System (11) quantifies the stress generated by the cells in response to the electrical stimulation, whose value at the collocation points is advanced in time using the same Runge–Kutta scheme adopted for the integration of the cell model.

Since the active stress depends on the potential, that value is computed at the collocation points via evaluation of the basis functions:

$$v\left(\tilde{\mathbf{r}}_{e_j}\right) = \mathbf{N}^e\left(\xi_j\right)\hat{\mathbf{v}} \quad j = 1, \ldots, n_e\,. \tag{85}$$

therefore coupling the monodomain formulation with the activation model. If the activation depends on the state variables $\mathbf{w}$ as well, these values must be calculated on the boundary of the domain to obtain the control variables.

In our approach, we compute the activation at every collocation point such that the corresponding control variables $\hat{\sigma}^{a_i}$ can be calculated by solving a square linear system:

$$\mathbf{N}^e\left(\xi_j\right)\hat{\sigma}^{a_i} = \sigma^a\left(\tilde{\mathbf{r}}_{e_j}, t_i\right) \quad j = 1, \ldots, n_e\,, \tag{86}$$

to define the field everywhere, such that the activation field can be evaluated in the shell cross-section corresponding to the quadrature point locations defined in the mechanical sub-problem. Having introduced our hypothesis, this operation results in a simple field evaluation:

$$\sigma^a(\theta^1, \theta^2, t_i) = \mathbf{N}^e\left(\theta^1_e, \theta^2_e\right)\hat{\sigma}^{a_i}\,, \tag{87}$$

even if the meshes are different since the parametrization of the domain is preserved under mesh refinement.

## 4. Numerical experiments

The numerical tests presented in this section explore and verify the features of the presented shell model. Initially, in Section 4.1, we analyze the quasi-static response of a cantilevered shell subjected to a prescribed active stress field. Afterward, in Section 4.2, we include the inertia forces, such that the time discretization is introduced as well. Finally, in Section 4.3, we propose a coupled electromechanical test, where the active stress is provided as a consequence of the instantaneous transmembrane potential field. In the latter test, we also include viscous damping to provide a complete description of the dynamic response. These tests mainly focus on tissue contractility to demonstrate the capabilities of the coupled numerical approach, therefore tissue pre-stretch is considered in a simplified model.

To simulate MTFs, we adopt a constitutive model calibrated on experimental data (Shim et al., 2012). Specifically, we consider a hyperelastic substrate made of polydimethylsiloxane (PDMS), modeled as an incompressible Neo-Hookean material, whose elastic energy reads:

$$\psi_n = \frac{1}{2}\mu\left(\text{tr}\left(\mathbf{C}\right) - 3\right)\,. \tag{88}$$

The second Piola–Kirchhoff stress tensor components, consequently, read (Kiendl et al., 2015):

$$S^{\alpha\beta}_n = \mu\left(\hat{g}^{\alpha\beta} - J^{-2}_0 g^{\alpha\beta}\right)\,, \tag{89}$$

with $J_0$ the in-plane Jacobian determinant:

$$J_0 = \sqrt{\frac{\det(g_{\alpha\beta})}{\det(\hat{g}_{\alpha\beta})}}. \tag{90}$$

Unlike previous works (Shim et al., 2012; Pezzuto et al., 2014; Böl et al., 2009), we assume that the material is completely incompressible rather than nearly incompressible. In fact, Kirchhoff Love shells are not prone to volumetric locking and do not require mixed formulations typical of 3D discretizations to prevent this kind of instability. For Neo-Hookean and other incompressible materials, the volumetric constraint can be enforced by elaborating analytically on the constitutive relation, as described in Kiendl et al. (2015).

In similarity with Shim et al. (2012), the passive response of the active layer is modeled using the same constitutive model enriched by an anisotropic component:

$$\psi_{ani} = \frac{\hat{E}_p}{\hat{a}^2}\left\{\exp\left[\hat{a}\left(\lambda_f - 1\right)\right] - \hat{a}\left(\lambda_f - 1\right) - 1\right\}\,, \tag{91}$$

being $\hat{E}_p$ and $\hat{a}$ the two material parameters characterizing the model while $\lambda_f$ the stretch in the fiber direction $\hat{\mathbf{f}}_0$. Accordingly, the anisotropic component of the second Piola–Kirchhoff stress tensor reads:

$$\mathbf{S}_{ani} = \frac{\hat{E}_p}{\hat{a}\lambda_f}\left\{\exp\left[\hat{a}\left(\lambda_f - 1\right)\right] - 1\right\}\left(\hat{\mathbf{f}}_0 \otimes \hat{\mathbf{f}}_0\right)\,, \tag{92}$$

or in curvilinear coordinates:

$$S^{\alpha\beta}_{ani} = \frac{\hat{E}_p}{\hat{a}\lambda_f}\left\{\exp\left[\hat{a}\left(\lambda_f - 1\right)\right] - 1\right\}\hat{f}^\alpha_0\hat{f}^\beta_0. \tag{93}$$

Summarizing, the second Piola–Kirchhoff stress tensor of the elastic substrate made of PDMS is given by the Neo Hookean model:

$$S^{\alpha\beta} = S^{\alpha\beta}_n\,, \tag{94}$$

**Table 1**
Parameters for the passive components of the constitutive laws.

| PDMS substrate: | | |
|---|---|---|
| $\rho_s$ | 0.965 | [mg/mm$^3$] |
| $\mu_s$ | 500 | [kPa] |
| Active layer: | | |
| $\rho_a$ | 0.965 | [mg/mm$^3$] |
| $\mu_a$ | 0.767 | [kPa] |
| $\bar{E}_p$ | 21 | [kPa] |
| $\bar{a}$ | 5.5 | [−] |

**Table 2**
Parameters for the IM approach.

| Mechanical activation model: | | |
|---|---|---|
| $\lambda_0$ | 1.24 | [−] |
| $\lambda_{min}$ | 0.86 | [−] |
| $\lambda_{max}$ | 1.34 | [−] |
| $\lambda_s$ | 1.14 | [−] |
| $\hat{P}$ | $\in [2.8, 21.6]$ | [kPa] |

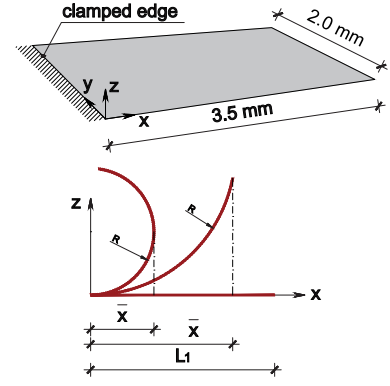

**Fig. 2.** Shell geometry and schematic representation of the procedure for curvature computation.

where the shear modulus $\mu_s$ is used, while, for the active layer, the Piola–Kirchoff stress tensor is given by the sum of several components:

$$S^{\alpha\beta} = S_n^{\alpha\beta} + S_{ani}^{\alpha\beta} + S_a^{\alpha\beta}. \qquad (95)$$

In this second material model, the shear modulus $\mu_a$ is used instead of $\mu_s$. We note that, in this work, we employ three integration points through the shell thickness in each material layer to compute normal forces and bending moments.

The active stress is calculated according to two different approaches, hereafter named *imposed* (IM) and *coupled* (CO) approaches. The first one neglects the electrophysiological stimulation but provides an activation pattern calibrated by means of experiments on *in-vitro* tissues. Conversely, the CO approach fully exploits the electromechanical coupling by means of an active stress field which locally depends on the action potential. However, it is calibrated for human cardiac tissues, not for *in-vitro* tissues. Nevertheless, the approximation is sufficient for numerical demonstration purposes. Therefore, we verify that realistic constitutive models, typical of the IM approach, are correctly reproduced, as well as complex implementations of coupled electromechanics, addressed following the CO approach.

In the IM activation model, the active stress $\sigma_a$ depends on the local stretch $\lambda_f$, varies in time according to a generic law $q(t)$ uniform in the domain, and the maximum value of stress is regulated by the parameter $\hat{P}$. In this work, we simplify the model proposed in Shim et al. (2012) including the pre-stretch of the active layer directly in the activation via a re-scaling of the relation between stretch and active stress, which results in:

$$\sigma_a = \begin{cases} \text{if } \lambda_{min} \leq \lambda_f \leq \lambda_{max}: \\ \hat{P} \, q(t) \left\{ 1 - \dfrac{\left[ \lambda_f + (\lambda_s - 1) - \lambda_0 \right]^2}{\left(1 - \lambda_0\right)^2} \right\} \\ \text{otherwise}: \\ 0 \end{cases}, \qquad (96)$$

being $\lambda_s$ the pre-stretch amplitude and $\lambda_0$ the optimal stretch, for which the active layer generates the maximum stress. Moreover, we remark that various formulations for $q(t)$ are used in the following examples to investigate different scenarios.

Conversely, in the context of a coupled electromechanical model (CO approach), the active stress depends on the transmembrane potential $v$ via the cardiac tissue model presented in Göktepe and Kuhl

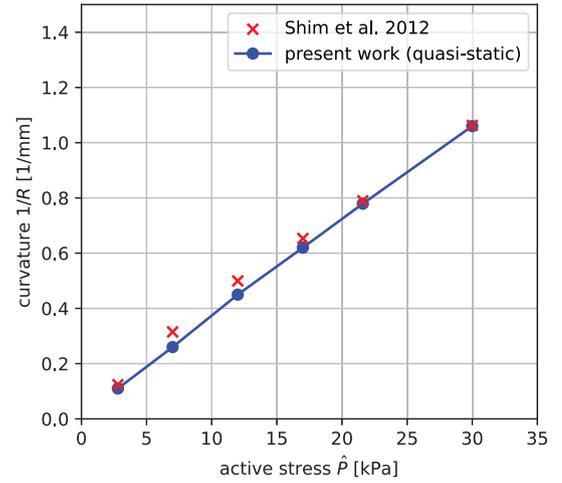

**Fig. 3.** Comparison between the quasi-static shell and the 3D model (Shim et al., 2012) in terms of curvature for several values of $\hat{P}$.

(2010). The evolution of the active stress is given by the following differential equation (cf. (11)):

$$\frac{\partial \sigma_a}{\partial t} = \zeta \left[ k_\sigma \left( v - v_r \right) - \sigma_a \right] \qquad (97)$$

being

$$\zeta = \zeta_0 + \left(\zeta_\infty - \zeta_0\right) \exp \left[ -\exp \left( -\xi_v (v - \bar{v}) \right) \right], \qquad (98)$$

coupled to the cell activity through the transmembrane potential $v$. The action potential $v$, in turn, depends on the tissue conductivity through (2) and, at the local level, on the ionic current (10). Herein, we model the electrophysiological activity via the normalized potential $v^*$, the time $t^*$, and the current $I^{ion*}$:

$$v^* = \frac{v - s_v}{r_v}, \quad t^* = \frac{t}{r_t}, \quad I^{ion*} = \frac{I^{ion} r_t}{r_v}, \qquad (99)$$

which appear in the following equations (cf. (10)):

$$\begin{cases} I^{ion*} = k_v \, v^* \left( v^* - a \right) - v^* w \\ \dfrac{\partial w}{\partial t^*} = \left( e_0 + \dfrac{\mu_{v1} w}{\mu_{v2} + v^*} \right) \left( -w - k_v v^* \left( v^* - b - 1 \right) \right) \end{cases}. \qquad (100)$$

For the CO approach, data are taken from the literature (Göktepe and Kuhl, 2010), while only the parameter $k_\sigma$ is re-scaled to achieve a maximum level of active stress equal to the average value of 12.2 kPa and the tissue conductivity is set to $D/C_m = 0.002 \, \text{mm}^2/\text{ms}$, to consider that the speed of propagation of the action potential in *in-vitro* tissues is reduced compared to *in-vivo* conditions. We re-scale by a factor of

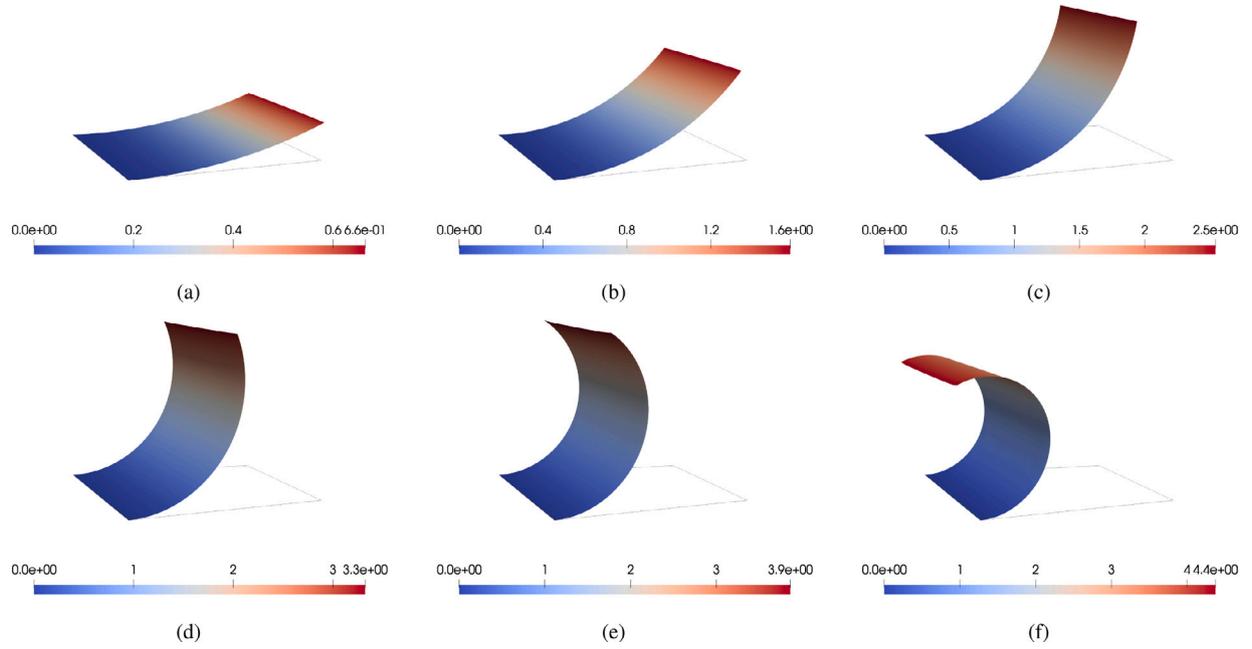

**Fig. 4.** Deflection of the cantilever analyzed in Section 4.1 for different values of maximum active stress. The displacement magnitude (in millimeters) is represented for: (a) $\hat{P}$ = 2.8 kPa, (b) $\hat{P}$ = 7 kPa, (c) $\hat{P}$ = 12 kPa, (d) $\hat{P}$ = 17 kPa, (e) $\hat{P}$ = 21.6 kPa, and (f) $\hat{P}$ = 30 kPa.

50 the value presented in Bueno-Orovio et al. (2008) to approximate a propagation speed in the order of 1 cm/s.

The complete list of parameters adopted in the present work is reported in Table 1 for the passive components of the constitutive model and in Tables 2–3 for the activation laws.

### 4.1. Quasi-static simulations with an imposed activation field

In this first example, we compare the proposed approach with 3D simulations from the literature (Shim et al., 2012). A thin shell ($L_1$ = 3.5 mm, $L_2$ = 2 mm, $d_s$ = 18 μm, $d_a$ = 4 μm), shown in Fig. 2, characterized by fibers oriented in the longitudinal direction $X$ ($\hat{f}_0 = (1,0,0)^T$), is uniformly activated following the IM approach in quasi-static conditions (i.e., the inertia terms are neglected) with different values of the maximum active stress $\hat{P}$ ($q(t) = 1$). The active layer induces a deflection of the cantilever quantified by the curvature ($1/R$), as defined in Shim et al. (2012), via the following relations:

$$\bar{x} = \begin{cases} R \sin\left(\dfrac{L_1}{R}\right) & \text{if } \bar{x} > \dfrac{L_1}{\pi} \\ R & \text{otherwise,} \end{cases} \tag{101}$$

being $\bar{x}$ the length of the projection of the deformed cantilever on the $XY$-plane computed at the mid-axis $Y$ = 1 mm, as shown in Fig. 2. Such a procedure was defined to compute the curvature starting from experimental measurements. For more details, the reader is referred to Shim et al. (2012). The structure is discretized with 50 × 10 quadratic elements and, to solve the non-linear problem via Newton's method, the activation is gradually imposed up to the maximum value in 42 steps.

Fig. 3 compares the results of the proposed shell-based approach to the active deformation of 3D simulations from the literature (Shim et al., 2012), highlighting a good matching, while Fig. 4 presents the corresponding deformed shapes of the shell mid-surface. As a remark, we note that similar results were achieved with a mesh refined twice.

In a second test, we fix the maximum value of active stress ($\hat{P}$ = 9 kPa) varying the thickness of the elastic substrate $d_s$ in the range 13÷28 μm. We compare the maximum curvature to the modified Stoney approximation (Shim et al., 2012):

$$\frac{1}{R} = \frac{6\left(3\mu_a + \hat{E}_p\right)\ln\lambda_0}{3\mu_s d_s}\left(\frac{d_a}{d_s}\right)\left(1 + \frac{d_a}{d_s}\right). \tag{102}$$

**Table 3**
Parameters for the CO approach.

| Mechanical activation model: | | |
|---|---|---|
| $k_\sigma$ | 0.122 | [kPa/mV] |
| $v_r$ | −80 | [mV] |
| $\zeta_0$ | 0.1 | [1/mV] |
| $\zeta_\infty$ | 1 | [1/mV] |
| $\xi_v$ | 1 | [1/mV] |
| $\hat{v}$ | 0 | [mV] |
| Electrophysiological model: | | |
| $k_v$ | 8 | [–] |
| $a$ | 0.15 | [–] |
| $b$ | 0.15 | [–] |
| $\mu_{v1}$ | 0.2 | [–] |
| $\mu_{v2}$ | 0.3 | [–] |
| $e_0$ | 0.002 | [–] |
| $r_v$ | 100 | [mV] |
| $s_v$ | −80 | [mV] |
| $r_f$ | 12.9 | [ms] |
| $D$ | 0.002 | [nF/ms] |
| $C_m$ | 1 | [nF/mm²] |
| Initial conditions at-rest: | | |
| $v_0$ | −80 | [mV] |
| $w_0$ | 0 | [–] |
| $\sigma_{a_0}$ | 0 | [kPa] |

Herein, we approximate the stretch in the active layer as uniform and equal to the model parameter $\lambda_0$, representing the optimal stretch for cardiomyocytes.

Numerical results in Fig. 5 are in agreement with the modified Stoney approximation: the same trend is observed, confirming the applicability of the proposed approach. However, we highlight that the proposed numerical approach has greater applicability since it relies on different hypotheses, for instance, material non-linearities and non-uniform activation can be considered. This last feature is fundamental when simulating complex activation patterns.

### 4.2. Dynamic simulations with an imposed activation field

The following tests aim at demonstrating the applicability of the methodology in an undamped dynamic framework.

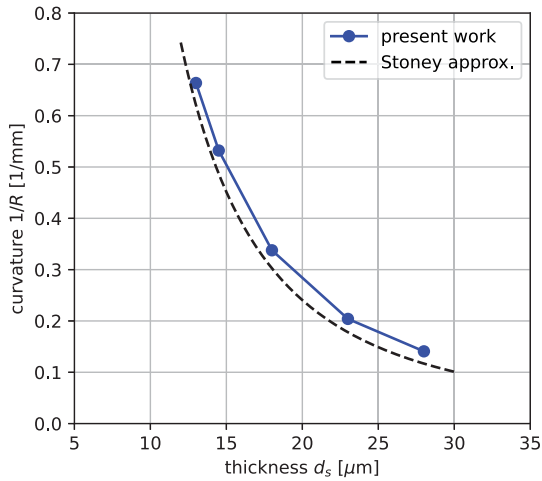

**Fig. 5.** Maximum curvature as a function of the substrate thickness for a fixed value of $\hat{P}$.

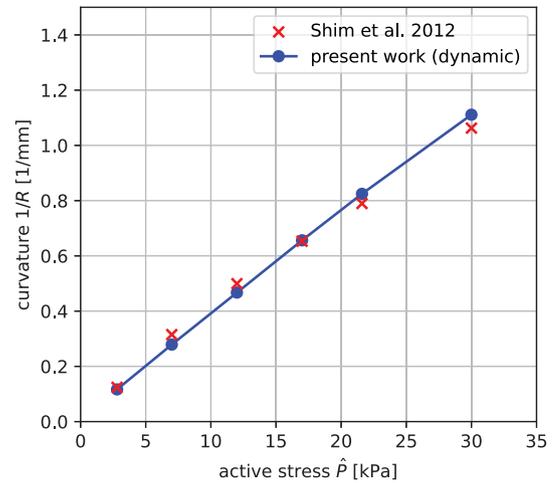

**Fig. 7.** Comparison between the dynamic shell and 3D model (Shim et al., 2012) in terms of curvature for several values of active stress $\hat{P}$, see also Fig. 3.

The cantilever simulated in Section 4.1 ($d_s = 18$ μm) is herein activated according to the timing-law proposed in (Shim et al., 2012):

$$q(t) = \left(\frac{t}{\bar{T}}\right)^2 \exp\left[1 - \left(\frac{t}{\bar{T}}\right)^2\right], \tag{103}$$

being $\bar{T} = 210$ ms, and characterized by a density $\rho$ equal to:

$$\rho = \rho_s \times d_s + \rho_a \times d_a, \tag{104}$$

where, for the sake of simplicity, the two densities are assumed equal to the density of PDMS. Starting from at-rest conditions (i.e., zero displacements and velocities), 500 ms of simulation are discretized using 100 generalized-$\alpha$ steps.

Displacements of the tip are shown in Fig. 6 for the maximum and minimum value of experimentally measured active stress ($\hat{P} = 2.8$ kPa and $\hat{P} = 21.6$ kPa) and compared to equivalent quasi-static simulations, for which the results at every instant are obtained through an independent analysis.

The results of the dynamic simulations are coherent with the quasi-static approach. Typical oscillations around the quasi-static results are observed and – as expected – a closer match is shown if the magnitude of the active stress is limited. Finally, we note that, in this dynamic example, the maximum curvatures are in agreement with the reference solution, as shown in Fig. 7, for all values of $\hat{P}$ analyzed in Section 4.1.

As a remark, it should be noted that in practical applications, the interaction of the MTF with the surrounding fluid (e.g., medium or imaging buffer) may introduce viscous dissipation. In the following tests, we include viscous effects by a damping matrix with constant coefficients to demonstrate the possibility of considering these factors in the simulations.

### 4.3. Coupled electromechanical simulations in the dynamic regime

In these examples, we demonstrate the possibility of effectively conducting coupled electromechanical simulations. Differently from previous tests, the activation field results from an electrophysiological simulation (CO approach) and, with such an aim, we adopt an available coupled-electromechanical activation model (Göktepe and Kuhl, 2010) in (97) instead of the a-priori imposed law (96).

Starting from at-rest conditions, see Fig. 8(a), the active layer is electrically stimulated (S1) at the clamped edge ($X = 0$ mm) by a distributed current flux ($I^n = 2$ mA/mm) for 5 ms. This generates a propagating action potential, which travels through the cantilever as shown in Fig. 8(b). At time $t = 430$ ms, a second stimulus (S2) is applied at the same site for 5 ms, as shown in Fig. 8(c). This second wave travels in a tissue that is not completely repolarized yet, see Figs. 8(d)–8(f), and, therefore, the action potential duration and the maximum force generated are reduced.

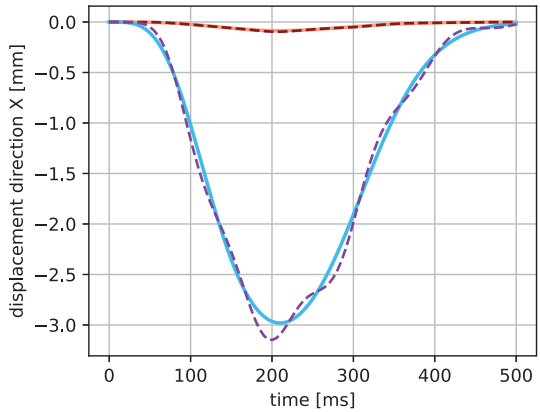

(a)

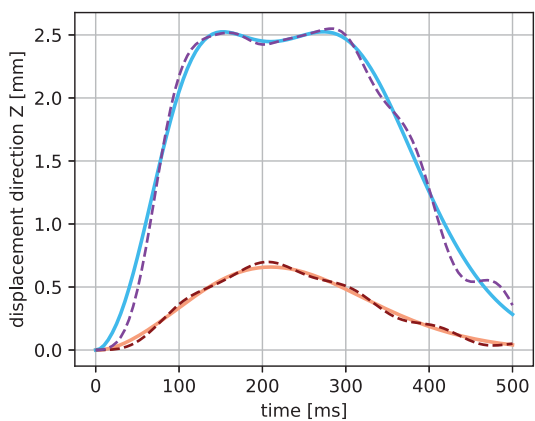

(b)

**Fig. 6.** Comparison between quasi-static and dynamic approaches. Displacements of a point located at $X = 3.5$ mm and $Y = 1$ mm (a) in the in-plane direction $X$, and (b) in the out-of-plane direction $Z$.

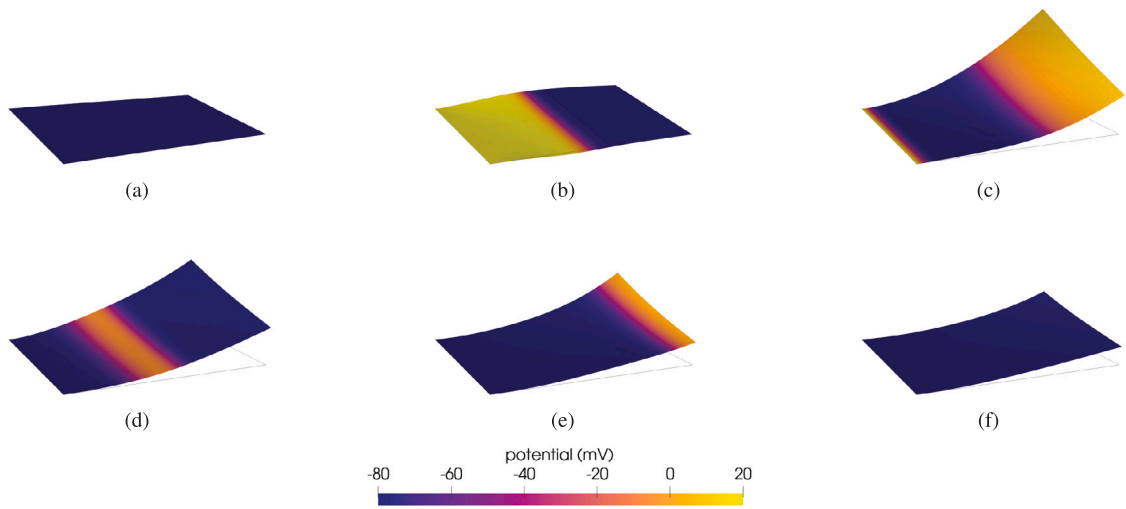

**Fig. 8.** Deflection of the cantilever undergoing the two planar waves at different times: (a) $t = 0$ ms, (b) $t = 100$ ms, (c) $t = 435$ ms, (d) $t = 575$ ms, (e) $t = 775$ ms, and (f) $t = 875$ ms. For the sake of simplicity, the transmembrane potential $v$ is plotted on the deformed shell mid-surface.

Given the activation field $\sigma_a(\hat{\mathbf{x}}, t)$ for $t \in [0, 2]$ s, the mechanical sub-problem is solved assuming that the fibers are oriented with a diagonal arrangement ($\hat{\mathbf{f}}_0 = (1, 1, 0)^T/\sqrt{2}$). Furthermore, as a proof of concept for viscous dissipation due to interactions with a fluid (Grosberg et al., 2011), we include in the dynamic formulation a damping matrix proportional to the mass:

$$\mathbf{C} = c\,\mathbf{M}, \tag{105}$$

with $c = 0.8$ m s$^{-1}$.

Different discretizations are used for each sub-problem: for the electrophysiological problem, we use $102 \times 60$ quadratic $C^1$ splines and $10^5$ time steps, while for the mechanics we use a coarser discretization composed by $50 \times 10$ quadratic elements and we subdivide the time interval in 400 steps.

In Fig. 9, we report the displacements and active stress $\sigma_a$ of the cantilever center (Point A: $X = 1.75$ mm, $Y = 1$ mm) and tip (Point B: $X = 3.5$ mm, $Y = 1$ mm). Results of a simulation conducted applying only the stimulus S1 are represented as well.

The coupled model automatically adapts to different stimulation conditions, varying the magnitude and the duration of the contraction according to the electrophysiological model; such a feature is further exploited in the last example, where the tissue undergoes spiral-wave excitations.

The stimulation protocol of this conclusive example is similar to the previous one: After the first planar wave, see Figs. 10(a)–10(b), a second stimulus (S3) is applied at time $t = 430$ ms, in the orthogonal direction (side $Y = 0$ mm) for 5 ms, as shown in Fig. 10(c). Such a combination of stimuli generates a self-sustained spiral wave, which rotates around a central region creating a complex activation pattern, as shown in Figs. 10(d)–10(h).

In Fig. 11, we report the displacements and active stress $\sigma_a$ of the cantilever center (Point A: $X = 1.75$ mm, $Y = 1$ mm), which is in the region around which the spiral-wave rotates, and tip (Point C: $X = 3.5$ mm, $Y = 2$ mm). Results of a simulation conducted applying only the stimulus S1 are represented as well.

In this last example, the activation pattern is even more complex and difficult to predict a-priori due to the autonomous firing activity triggered by the spiral wave. This test confirms that coupled simulations are necessary in scenarios where different stimuli interact.

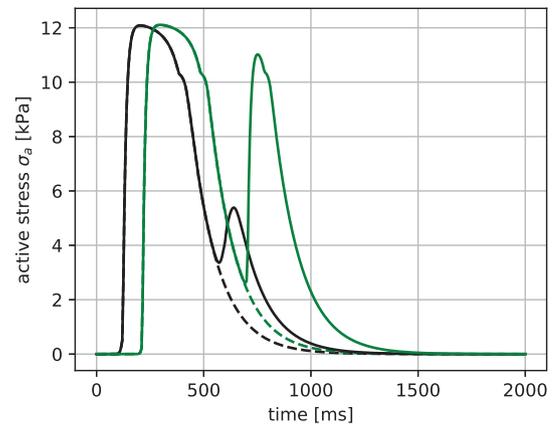

(a)

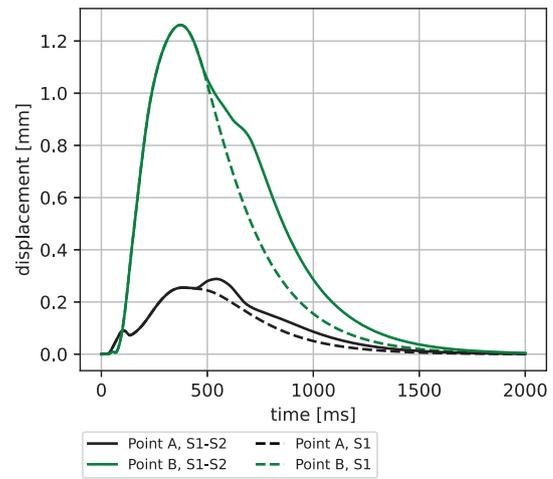

(b)

**Fig. 9.** Active stress (a) and corresponding displacement magnitude (b) in two points (A: center, B: tip) of the cantilever subject to two planar-waves (S1–S2) and comparison with results of a single planar-wave (S1).

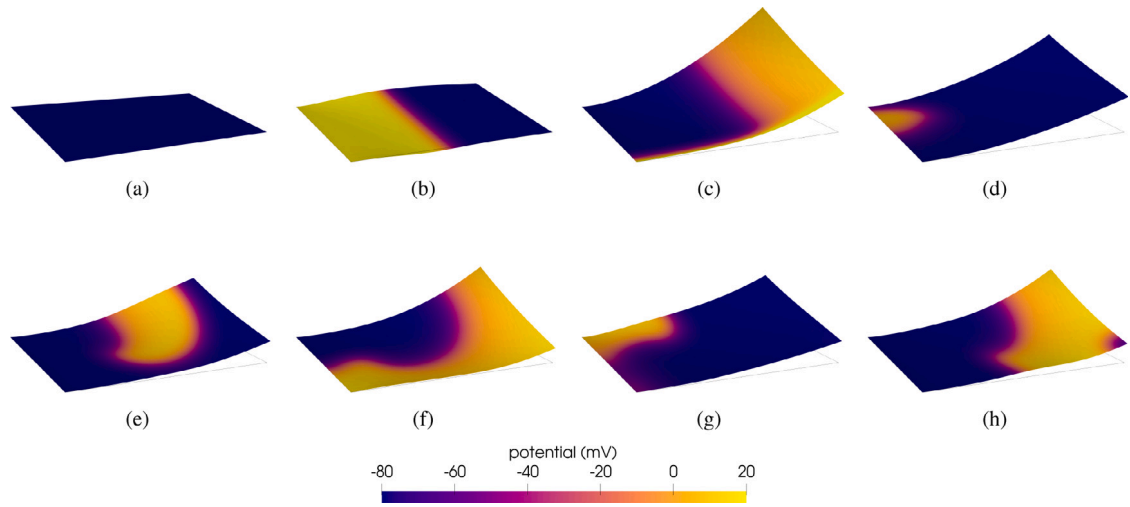

**Fig. 10.** Deflection of the cantilever undergoing the spiral-wave at different times: (a) $t = 0$ ms, (b) $t = 100$ ms, (c) $t = 435$ ms, (d) $t = 575$ ms, (e) $t = 775$ ms, (f) $t = 875$ ms, (g) $t = 1010$ ms, and (h) $t = 2000$ ms. For the sake of simplicity, the transmembrane potential $v$ is plotted on the deformed shell mid-surface.

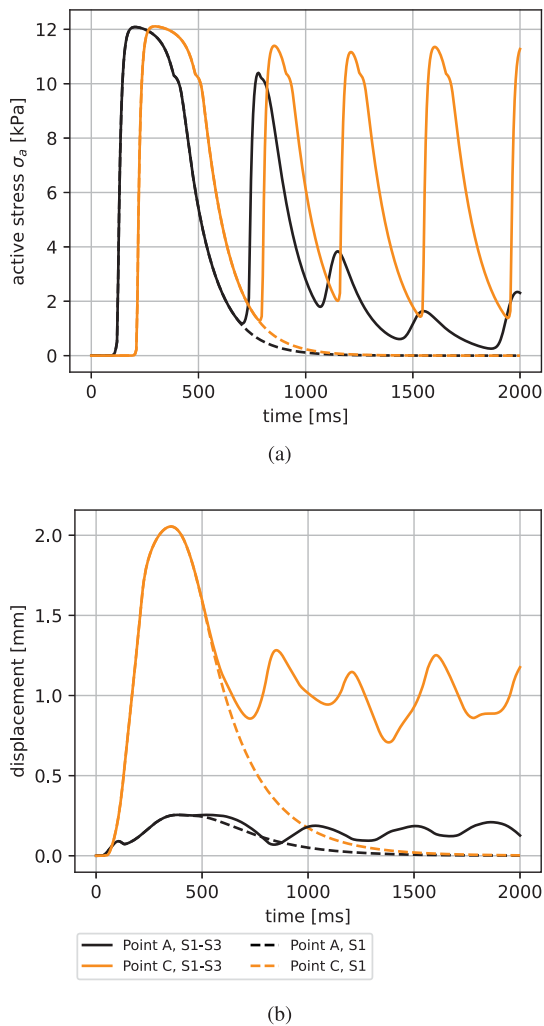

**Fig. 11.** Active stress (a) and corresponding displacement magnitude (b) in two points (A: center, C: corner) of the cantilever undergoing the spiral-wave (S1–S3) and comparison with results of a single planar-wave (S1).

## 5. Conclusions and outlooks

In the present work, we have proposed an active stress approach to model the activation of artificial MTFs relying on thin shells instead of full 3D models. Specifically, we have presented how to compute the internal forces and the bending moments for these kinds of structures where the materials differ in the passive and active properties relying on numerical integration through the thickness, which is needed in any case if the materials exhibit nonlinear behaviors. Numerical tests have proven that the proposed approach is suitable for both quasi-static and dynamic conditions.

The cell activation may derive from an electrophysiological simulation, where the cell activity is modeled by solving the monodomain formulation on a bi-variate manifold representing the contractile layer. Moreover, the electromechanical coupling has been addressed. We have demonstrated that the classical activation model can be re-used in this context without modifications and that the electrophysiological problem can be discretized using different approaches. Indeed, the Isogeometric-Galerkin method is used to discretize the mechanical field, while the Isogeometric-Collocation method is used for the electrical part of the problem. Moreover, different meshes and time step sizes can be used for each sub-problem, further demonstrating the possibility of tailoring the discretization to enhance the numerical performance.

Compared to a full 3D model, the proposed approach reduces the number of degrees of freedom required to conduct the simulations since only the displacements of the midsurface are explicitly computed. A similar rationale holds for the electrophysiological problem. Moreover, the isogeometric approach easily enables high-order discretizations, which can achieve high accuracy for a limited number of degrees of freedom. Finally, the use of a collocation scheme avoids quadrature in favor of evaluations at collocation points (whose number is equal to that of control points), significantly limiting the effort in the integration of the cell model.

Currently, the mesh coupling relies on a global projection. Future studies may improve this part of the simulation process by exploiting local quasi-interpolation operators (Buffa et al., 2016) to parallelize the process. Furthermore, the proposed approach can be recast in the immersed framework (Torre et al., 2023b; Loibl et al., 2023) to decouple the topology of parametric space and physical domain.

The proposed methodology, for instance, can be employed to evaluate the state of stress at the interface between the active layer and substrate to evaluate the risk of unbinding. Indeed, the computed displacement field can be post-processed via stress recovery procedure (Patton et al., 2021b,a, 2019) to evaluate both the in-plane and out-of-plane components of the Cauchy stress tensor. Such values can be successively compared with the maximum resistance of the connection. Eventually, from an application point of view, the approach can be used to design actuators and power generation devices (Feinberg et al., 2007) or to optimize the shape of soft robots (Lucantonio et al., 2014; Ricotti and Fujie, 2017).

**Declaration of competing interest**

The authors declare the following financial interests/personal relationships which may be considered as potential competing interests: Alessandro Reali reports financial support was provided by CENTRO NAZIONALE HPC, BIG DATA E QUANTUM COMPUTING. Michele Torre reports financial support was provided by CENTRO NAZIONALE HPC, BIG DATA E QUANTUM COMPUTING. Alessandro Reali reports financial support was provided by Italian Ministry of University and Research. Michele Torre reports financial support was provided by Italian Ministry of University and Research. Francesco S. Pasqualini reports financial support was provided by European Research Council. Josef Kiendl reports financial support was provided by European Research Council. Francesco S. Pasqualini reports financial support was provided by National Center for Gene Therapy and Drugs based on RNA Technology. Marco D. de Tullio reports financial support was provided by Italian Ministry of University and Research. If there are other authors, they declare that they have no known competing financial interests or personal relationships that could have appeared to influence the work reported in this paper.

**Acknowledgments**

AR and MT acknowledge the contribution of the National Recovery and Resilience Plan, Mission 4 Component 2 Investment 1.4 CN_00000013 CENTRO NAZIONALE HPC, BIG DATA E QUANTUM COMPUTING, spoke 6. The support of the Italian Ministry of University and Research (MUR) through the PRIN project COSMIC (No. 2022A79M75), funded by the European Union - Next Generation EU, is also acknowledged. FSP is partially supported by the European Research Council (ERC) through the starting grant SYNBIO.ECM (no. 852560) and by the National Recovery and Resilience Plan, Mission 4 Component 2 Investment 1.4 CN_00000041 Metabolic and cardiovascular diseases, spoke 4. JK is partially supported by the European Research Council through the H2020 ERC Consolidator Grant 2019 n. 864482 FDM². MDdT has been partially supported by MUR PRIN 2022 project ABYSS: Accurate simulation of Bio-hYbrid Soft Swimmers, Grant No. 2022BN3RAC, CUP: D53D23003410006, funded by European Union - Next Generation EU.

**Appendix A. Notes on curvilinear coordinates**

The coordinates of a generic point $\mathbf{x}$ can be expressed in cartesian coordinates as:

$$\mathbf{x} = x^i \mathbf{e}_i \quad i = 1, 2, 3 \,, \tag{A.1}$$

where $\mathbf{e}_i$ is the $i$th cartesian base vector, or equivalently as:

$$\mathbf{x} = \theta^i \mathbf{g}_i = \theta_i \mathbf{g}^i \quad i = 1, 2, 3 \,, \tag{A.2}$$

$\mathbf{g}_i$ are the covariant base vector and $\mathbf{g}^i$ are the contravariant base vector.

The covariant vectors are defined as:

$$\mathbf{g}_i = \frac{\partial \mathbf{x}}{\partial \theta^i} \,, \tag{A.3}$$

while the contravariant vectors are such that the following relation holds:

$$\mathbf{g}_i \cdot \mathbf{g}^j = \delta_i^j \,, \tag{A.4}$$

being $\delta_i^j$ the Kronecker Delta symbol.

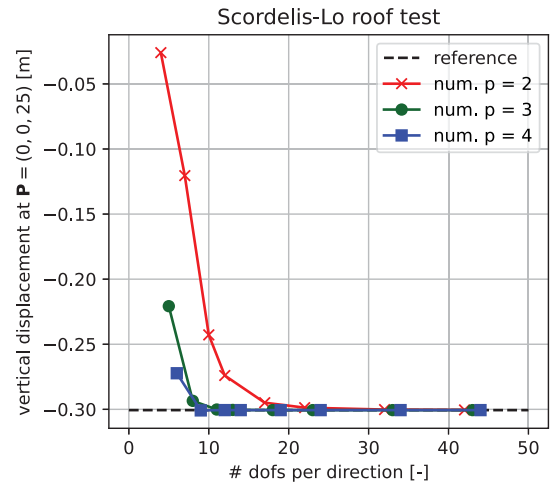

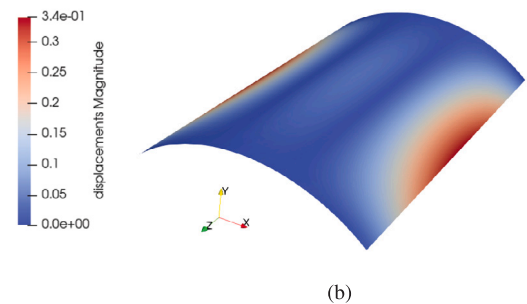

**Fig. B.1.** Results of the Scordelis–Lo roof test. (a) Convergence analysis. (b) Displacement field simulated using 40 quartic $C^3$ elements per direction.

In relating the two bases, one of the key quantities is the metric tensor $\mathbf{g}$, which can be expressed as:

$$\mathbf{g} = g^{ij} \mathbf{g}_i \otimes \mathbf{g}_j = g_{ij} \mathbf{g}^i \otimes \mathbf{g}^j \,, \tag{A.5}$$

being the coefficients related by the matrix inversion operation:

$$[g^{ij}] = [g_{ij}]^{-1} \tag{A.6}$$

and

$$g_{ij} = \mathbf{g}_i \cdot \mathbf{g}_j \,. \tag{A.7}$$

The previous expressions provide relations between the covariant and contravariant vectors:

$$\mathbf{g}^i = g^{ij} \mathbf{g}_j \,, \tag{A.8}$$

$$\mathbf{g}_i = g_{ij} \mathbf{g}^j \,. \tag{A.9}$$

Furthermore, derivatives of contravariant vectors with respect to the coordinates $\theta^i$ are computed, according to Itskov et al. (2007), by exploiting the Christoffel symbol of second kind $\Gamma_{jk}^i$:

$$\mathbf{g}_{,k}^i = -\Gamma_{jk}^i \mathbf{g}^j \,, \tag{A.10}$$

being

$$\Gamma_{jk}^i = \mathbf{g}_{j,k} \cdot \mathbf{g}^i \,. \tag{A.11}$$

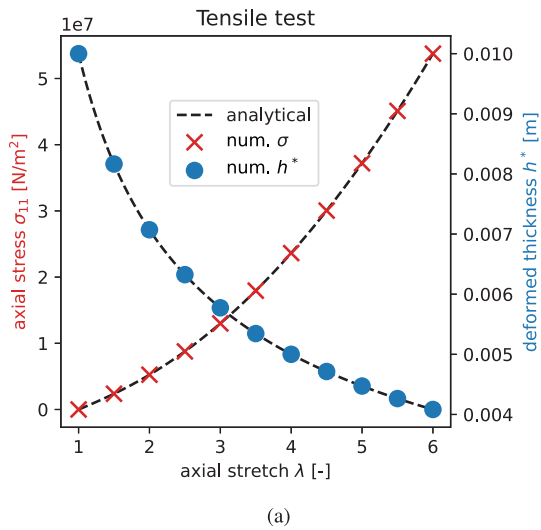

(a)

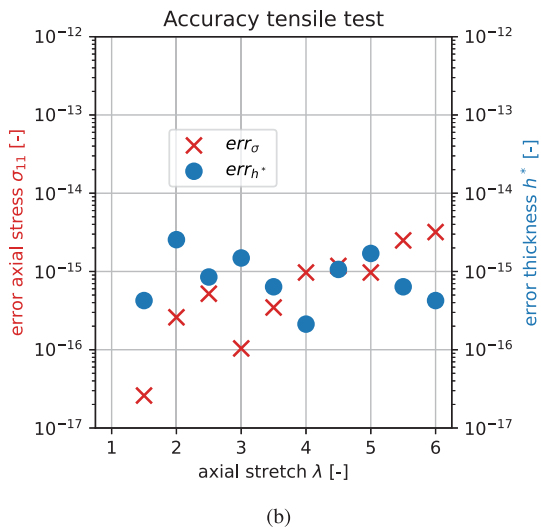

(b)

**Fig. B.2.** Nonlinear tensile test: (a) magnitude of the axial stress and shell thickness. (b) Relative error for the two quantities with respect to the analytical solution.

## Appendix B. Verification of the structural solver

We verify the implementation of the Kirchhoff–Love shell discretization by solving benchmarks and comparing the results against analytical solutions.

First, we perform linear tests by reproducing the Scordelis–Lo roof test, considering the data reported in Kiendl et al. (2009), for which the analytical solution of the Y-displacement in the middle of the straight free edge is provided. In Fig. B.1, results of the convergence test under mesh refinement for B-splines of degrees two, three, and four are shown.

Successively, the non-linear regime is tested. We reproduce the traction test presented in Kiendl et al. (2015) using an incompressible Neo-Hookean material. The results of stress and contraction in the normal direction in the middle of the domain are reported. The results of a simulation performed using four quadratic $C^1$ elements are in agreement with the analytical solution, as demonstrated by the magnitude of the relative error computed in the middle of the domain, shown in Fig. B.2(b). Furthermore, we note that the magnitude of other components of the stress tensor, which should be null, are in the order of $10^{-9}$ N/m$^2$.

Both tests confirm the correctness of the implementation.

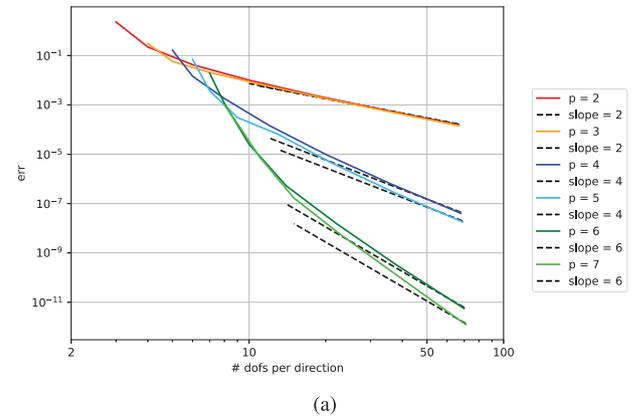

(a)

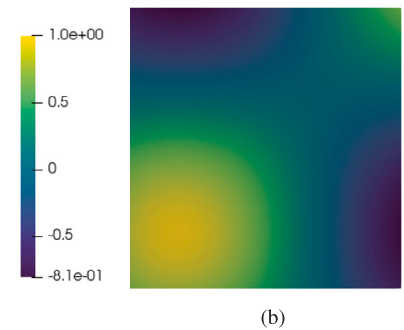

(b)

**Fig. C.1.** Results of the Poisson problem. (a) Convergence analysis. (b) Field simulated using 64 knot spans per direction, $p = 7$, and maximum continuity.

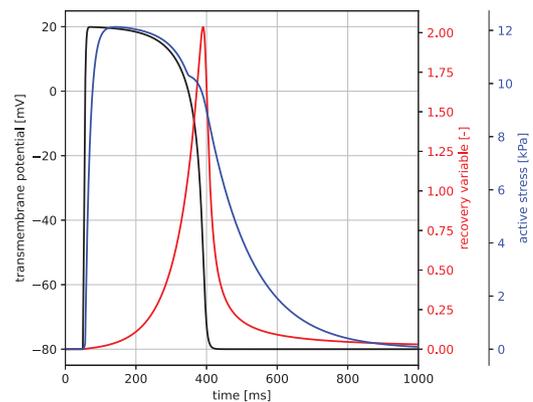

**Fig. C.2.** Results of a single cell simulation.

## Appendix C. Verification of the electrophysiological solver

We verify the correctness of the collocation approach by checking the rate of convergence of the $L^2$-norm of the error for a problem where the analytical solution is available. To this end, we focus on the Poisson equation:

$$\Delta v = q, \tag{C.1}$$

considering it as a proxy of the reaction–diffusion partial differential equation.

**Table C.1**
Control point coordinates defining the geometry of the domain via a single bi-quadratic knot span.

| X | Y | Z |
|---|---|---|
| 0.3 | 0.3 | 0 |
| 0.8 | 0.3 | 0 |
| 1.3 | 0.3 | 0 |
| 0.3 | 0.8 | 0 |
| 0.9 | 0.9 | 0 |
| 1.3 | 0.8 | 0 |
| 0.3 | 1.3 | 0 |
| 0.8 | 1.3 | 0 |
| 1.3 | 1.3 | 0 |

Applying Neumann boundary conditions on the right, top, and bottom sides and Dirichlet boundary conditions on the left edge of the domain in Fig. C.1, we conduct convergence studies exploiting the following definition of the error:

$$err = \sqrt{\frac{\int_{A_e} (v - v^a)^2 \, dA}{\int_{A_e} (v^a)^2 \, dA}}, \tag{C.2}$$

being $v^a$ the analytical solution

$$v^a = sin(\pi x) \times sin(\pi y). \tag{C.3}$$

To better verify the implementation, we ensure that the covariant vectors are not orthogonal by deforming the mesh. Coordinates of the control points used to represent the geometry are presented in Table C.1.

Under mesh refinement, we obtain the expected rates of error convergence for all the degrees analyzed.

Finally, we perform a simulation of the single-cell activity, shown in Fig. C.2, to verify the implementation of the cell and activation models. Starting from at-rest conditions, the cell is stimulated by an applied current for 2 ms after 50 ms of activity to trigger an action potential and the corresponding mechanical activation, whose evolution is simulated for 1 s using the same scheme adopted in the previous numerical examples. Results are in agreement with similar simulations presented in a previous work (Göktepe and Kuhl, 2010).